\documentclass{amsart}
\usepackage[T1]{fontenc}
\usepackage[cp1250]{inputenc}
\usepackage[centertags]{amsmath}
\usepackage{amsfonts}
\usepackage{amsthm}
\usepackage{newlfont}
\usepackage{amscd}
\usepackage{amsgen}
\usepackage{pb-diagram}
\newlength{\defbaselineskip} 
\theoremstyle{plain}
\newtheorem{thm}{Theorem}[section]
\newtheorem{prop}[thm]{Proposition}
\newtheorem{cor}[thm]{Corollary}
\newtheorem{conj}[thm]{Conjecture}
\newtheorem{lem}[thm]{Lemma}

\theoremstyle{remark} \newtheorem{rem}[thm]{Remark}
\theoremstyle{definition} 
\theoremstyle{definition} \newtheorem{ex}[thm]{Example} %

 \numberwithin{equation}{section}
\newcommand{\op}{\operatorname}

\newcommand{\ra}{\rightarrow}

\begin{document}

\title{Correspondences between modular Calabi--Yau fiber products}
\author{Micha{\l} Kapustka}
\thanks{The project is co-financed from the European Union funds and national budget.}
\maketitle
\begin{abstract}
We describe two ways to construct finite rational morphisms
between fiber products of rational elliptic surfaces with section
and some Calabi--Yau manifolds. We use them to construct
correspondences between such fiber products that admit at most
five singular fibers and rigid Calabi--Yau threefolds.
\end{abstract}

\section*{Introduction}
It is well known that if there exist a correspondence defined over
$\mathbb{Q}$ between two varieties then the L-series of these
varieties admit a common component. In fact the Tate conjecture
states that if two varieties have L-series admitting a common
component, then there should exist a correspondence between them.
It is also a well known fact that almost all rigid Calabi--Yau
threefold are modular (see \cite{DieuMaroh,Dieu2Mar}). On the
other hand modularity for nonrigid Calabi--Yau threefolds is still
a widely open question. Most of the known examples of modular
nonrigid Calabi--Yau threefolds (see
\cite{Liv3Yui,HV,Hulek,Schuett,Schuett2}) are in some way
connected with fiber products of elliptic surfaces with section.
For instance Hulek and Verrill (in \cite{Hulek}) gave a method to
prove the modularity of some nonrigid Calabi--Yau manifolds based
on the existence of many elliptic ruled surfaces. The authors
prove the modularity of a large class of Calabi--Yau fiber
products with the so-called modular defect equal to 0 (i.e. being
a fiber product of non isogenous surfaces that admits exactly five
singular fibers). They also find numerical evidence for the
newforms of weight 4 in many of the examples.

In this paper we investigate some constructions of correspondences
for Calabi--Yau threefolds which are resolutions of fiber products
of rational elliptic surfaces with section. These constructions
give us new tools to study the Tate conjecture for the concerned
varieties. In particular they permit us to find correspondences
between varieties that are conjectured by Hulek and Verrill to
have the same modular forms.

Moreover for most explicit examples it is not a hard task to check
if constructed correspondences are defined over $\mathbb{Q}$. In
this way we get a new proof of modularity for these examples based
on the theorem of Dieulefait and Manoharmayum. In particular we
can check modularity for all examples described in the tables from
\cite{Hulek}. In many cases our results together with the method
of Faltings--Serre--Livn\'e are also useful to prove that the
predicted newform of weight $4$ is indeed correct (for a specific
example see \cite{AnnMG}).

We study two types of correspondences. The first type consist of
correspondences between a fiber products of rational elliptic
surfaces with section and threefolds arising from these products
by a fiberwise Kummer construction. This construction is well
known and described in \cite{Schoen}. What we only need to add is,
when such Kummer constructions lead to rigid varieties. To find
these we use results of \cite{CV} and \cite{MG} describing
deformations of such Kummer fibrations by interpreting them as
double octics. The second construction is based on the fact that
an isogeny between two elliptic surfaces $S_1$ and $S_2$ induces
an isogeny between the fiber products of these surfaces with any
other elliptic surface with section. To obtain these
correspondences we study isogenies between elliptic surfaces.
\section*{Acknowledgements}
I would like to thank my advisor S. Cynk for many fruitful
discussions during the work on this subject.
\section{Isogenies between elliptic surfaces.}
We start with the study of isogenies between relatively minimal
rational elliptic surfaces with section over $\mathbb{P}^1$. By
the term isogeny we mean an isogeny of these surfaces seen as
elliptic curves over the function field of $\mathbb{P}^1$. Let $S$
denote such a rational elliptic surface with chosen 0 section. Let
$\pi:S\ra \mathbb{P}^1$ be its canonical projection. We know that
the set of smooth points $F^{\sharp}$ of a fiber $F$ has a natural
structure of group induced by the structure of the group of
sections.

\subsection{Isogenies of order 2.} There are two types of birational involutions preserving fibers
acting on the surface $S$.
\begin{itemize}
 \item The first involution acts on the set of smooth points of any fiber by $x\mapsto b(\pi(x))-x$, where $b$ is a given
section
 \item The second involution acts on the set of smooth points of any fiber by $x\mapsto x+a(\pi(x))$, where $a$
is a given nontrivial section such that $a+a=0$.
\end{itemize}
Assume moreover that $S$ has only reducible fibers. Then the
indeterminacy loci for both types of involutions are finite sets
of points. Thus these birational involutions extend to
automorphisms of the surfaces. The first type of involution is
used in the Kummer construction on the fiber product. It's fixed
locus is a four-section on $S$ given by $\{x \colon
x+x=b(\pi(x))\}$. Whereas an involution of the second type has no
fixed points on the set of smooth points of the fibers. Hence the
resolution of the quotient of $S$ by such an involution gives
another elliptic surface with section. Unfortunately the
two-torsion sections $a$ do not always exist. In \cite{Miranda}
the author gives two sufficient conditions for the existence of a
torsion section in terms of types of singular fibers.
\begin{lem}[\cite{Miranda}] \label{warunek na istnienie sekcji torsyjnej}
Let us assume that $S$ has $n$ semi-stable fibers.  Let us denote
the indices of the singular fibers by $k_1,\dots, k_n$. Let $p$ be
a prime number. Let one of the following cases holds:
\begin{enumerate}
 \item The prime $p$ divides the indices $k_i$ for $i\geq 4$,
 \item The prime $p$ divides the indices $k_i$ for $i\geq 5$ and
 moreover one of the following cases holds
 \begin{enumerate}
 \item The prime $p=2$. The numbers $k_i$ are divisible by $4$ for
$i\geq5$. The following equality holds
$$(-1)^n k_1k_2k_3k_4 \not= \prod_{i\geq5} (k_i -1) \quad
\operatorname{mod}8.$$
 \item The prime $p$ is odd and $k_1k_2k_3k_4$ is not a square
 modulo $p$.
 \end{enumerate}
\end{enumerate}
Then there is a $p$ torsion section on $S$.
\end{lem}

From a more pictorial point of view we can see the given surface
$S$ as the resolution of a double cover of $\mathbb{F}_2$ ramified
over the sum of the section and a three-section $T$ disjoint form
it (see \cite{MirPer}). We see that a section $a$ to satisfy
$a+a=0$ and $a\not\equiv0$ has to be the proper transform (by the
resolution of singularities) of a section contained in the divisor
$T$. Such a section exists if and only if $T$ is reducible.
Observe moreover that if $S$ is an elliptic surface with a nonzero
section $a$ such that $a+a=0$, then the quotient surface $S'$ of
$S$ by the involution $x\mapsto x+a$ also admits such a section.
Indeed the sections $a$ and $0$ are mapped into the zero section
of $S'$ and the image of the remaining two-section has to be a
section $d$ on $S'$ satisfying $d+d=0$. Finally, if the surface
$S$ is rational, then by the L\"uroth theorem $S'$ also is.

We have just proven the following lemma.
\begin{lem} Let $S$ be a rational elliptic surface with chosen zero
section. The following conditions are equivalent:
\begin{itemize}
 \item There exists a 2:1 isogeny $f$ between $S$ and
another rational elliptic surface $S'$ with section such that $f$
maps fibers onto fibers.
 \item The surface $S$ admits a two-torsion section.
 \item The fixed locus of the involution $x\mapsto -x$ has at least
 three components (it has always at least 2).
\end{itemize}
\end{lem}
We can see that the number of components of the branch curve
cannot always be deduced from the types of singular fibers. In
fact we can find two surfaces with the same types fibers, the
first admitting an isogeny while the second not.
\begin{ex} Let $Q_1$ be an irreducible quartic curve in $\mathbb{P}^2$
with three nodes. Let $P$ be a generic smooth point on $Q_1$. Let
$\hat{S}_1$ be the double cover of $\mathbb{P}^2$ branched over
$Q_1$. Let $S_1$ be the resolution of the singularities of the
blow up of the surface $\hat{S}_1$ in the inverse image of $P$.
Then $S_1$ is an elliptic surface with singular fibers
$I_2I_2I_2I_2I_1I_1I_1I_1$. If we take the inverse image of $P$ as
the zero section then the corresponding three-section $T$ on
$\mathbb{F}_2$ will be birational to the quartic $Q_1$, hence
irreducible. If we perform the same construction for $Q_2$ being
the sum of a smooth cubic and a line cutting it transversely we
will obtain a surface with the same type of fibers, but inducing a
reducible divisor $T$.
\end{ex}
The following lemma gives us a necessary condition for the surface
$S$ to admit a torsion section.
\begin{lem}\label{moje war na istnienie sekcji polowkowej}
Let us assume $S$ has only semi-stable fibers. Let us denote by
$k_i$ the indices of the singular fibers. If the number of odd
indices $k_i$ is $>4$ then the surface does not admit any
two-torsion section.
\end{lem}
\begin{proof} Recall that the indices $k_i$ of singular fibers correspond to suitable singularities of the curve $T$.
Observe moreover that the assumption of the lemma is equivalent to
the inequality $\sum \lfloor\frac{k_i}{2}\rfloor<4$. The proof of
is now straightforward by contraposition. Assume that $T$ is
reducible then the index of intersection of two components of $T$
is 4. Indeed, let $T=a+d$. We compute the intersection index $a
\cdot b$ by writing $a$ and $d$ in terms of the zero section $s$
and the fiber $F$. As $a$ is a section disjoint from $s$ we have
$a=s+2F$. Similarly $d=2s+4f$. This gives $a \cdot d= 4$. This
translates to indices of singularities of the curve $T$ and hence
to the $k_i$ giving a contradiction.
\end{proof}

Let us now concentrate on the case where the isogeny does exist.
We keep the notation from the beginning of the section. We are
interested in studying the image of the isogeny. To do this we
need to know the action of the involution $x\mapsto x+a$ on fibers
of $S$.

\begin{lem}
Assume that $S$ is an elliptic surface with section with reduced
fibers and admitting an isogeny of order 2. Let $i$ be the
involution associated to the isogeny. Let $F$ be a fiber of $S$.
We have the following possibilities:
\begin{enumerate}
 \item The fiber $F$ is smooth. Then $F/i$ is also a smooth elliptic curve.
 \item The fiber $F$ is of type $I_{2k+1}$. Then $i$ acts on each irreducible
 component of $I_{2k+1}$ separately and keeps fixed their
 intersection points.
 \item The fiber $F$ is of type $I_{2k}$. Then we have two possibilities:
 \begin{enumerate}
 \item The involution $i$ acts on each irreducible
 component of $I_{2k}$ separately and keeps fixed their
 intersection points.
 \item The involution $i$ interchanges pairs of opposite
 irreducible components and has no fixed points.
 \end{enumerate}

 \item The fiber $F$ is of type $III$. Then the involution
interchanges the two components of $F$ and fixes their
intersection point.
\end{enumerate}
In particular $F$ cannot be of type $IV$ or $II$
\end{lem}
\begin{proof}

Recall that to each fiber of $S$ there is an associated groups
$F^{\sharp}$ of smooth points of  the fibers. Moreover there is
also the group $F_0$ of smooth points of the component cutting the
zero section. This gives rise to another group $F^{\sharp}/F_0$
the group of components of the fiber. The table describing these
groups is given in \cite[Table 4]{Shaf2}. The group
$F^{\sharp}/F_0$ is cyclic for all reduced fibers $F$ whereas the
group $F_0$ is either the multiplicative group $\mathbb{C}^*$ or
the additive group $\mathbb{C}$ depending on whether the fiber $F$
is semi-stable or not. Using this we observe that there never
exists a section $a$ satisfying our assumptions if $S$ admits a
fiber of type $IV$ or $II$. This follows from the fact that
neither $\mathbb{Z}/3\mathbb{Z}$ nor the additive group
$\mathbb{C}$ admit nonzero fixed points of the inversion. This
also proves that the only possibility for the involution on the
fiber of type $III$ is to act nontrivially on the set of
components of the fibers. This ends the case of non semi-stable
fibers. Let $F$ be now of type $I_n$.

For $n$ odd the group $\mathbb{Z}/n\mathbb{Z}$ does not admit a
nontrivial element of order 2, hence the action of the involution
on this group has to be trivial. Moreover, a translation on the
multiplicative group $\mathbb{C}^*$ has no fixed points hence the
points of intersection of the component have to remain fixed.

For $n$ even we have two possibilities.
\begin{itemize}
\item The section $a$ cuts $F^0$ the component of $F$ cutting the
zero section. Then $i$ acts on each component separately and
analogically to the case of $n$ odd.

\item The section $a$ does not cut the component $F^0$. Then the
involution $i$ interchanges opposite components. It obviously does
not have any fixed points for $n>2$. For $n=2$ we use the fact
that the quotient $F/i$ does also have to be a special fiber of an
elliptic fibration.
\end{itemize}
\end{proof}

Let us consider the surface $S'$ which is the minimal resolution
of the image of $S$ by the isogeny. We can see $S'$ as the minimal
resolution of the quotient $S/i$ of the surface $S$ by the
involution $i$.

Observe that the fixed points of the involution induce nodes on
the quotient variety. After resolving these singularities by
blowing them up the fibers will have more components. More
precisely in this case fibers of type $I_n$ on $S$ will correspond
to fibers of type $I_{2n}$ on $S'$.

We have just seen that the studied involution $i$ changes only the
type of semi-stable fibers. That is, the surface $S'$ which is the
minimal resolution of the quotient $S/i$ has singular fibers in
the same place as the singular fibers of $S$, but at the place of
a fiber $I_n$ there is a fiber $I_{2n}$ or $I_{n/2}$. Moreover, we
know that a fiber of $S'$ is $I_{2n}$ if and only if the section
$a$ meets the component $F^0$ of the fiber $F=I_n$. In the terms
of the trisection $T=a + d$ it means that $F$ does not pass
through the intersection point of the section $a$ with the
two-section $d$. In terms of the singular fibers of $S$ this means
that there are $\leq 4$ semi-stable fibers whose indices are even
and sum up to at most eight that are mapped to fibers with indices
divided by two. The indices of the remaining semi-stable fibers
are multiplied by 2, and the fibers of type III map to fibers of
type III.

\subsection{Higher order isogenies on elliptic surfaces.}
Observe that most of the above can also be formulated in the case
of higher order isogenies. More precisely each section with
torsion of degree $p$ for $p$ prime induces an action of a cyclic
group on $S$. A sufficient condition for the existence of such a
section for a fixed $p$ is given in Lemma \ref{warunek na
istnienie sekcji torsyjnej}. We can also analogously prove the
following lemma.
\begin{lem} We have two possibilities for the action of the cyclic
group of rank $p$ onto a semi-stable fiber $F$ of an elliptic
surface $S$.
\begin{enumerate}
\item The action has no degenerate orbits. Then $F$ is of type
$I_{p.k}$ for some number $k$ and the induced action on the group
$F^\sharp$ is cyclic.

\item The action is trivial on all points of intersection of the
components of the fiber.
\end{enumerate}
\end{lem}

Observe that in the second case the quotient of the surface $S$ by
the group action admits singularities in the fixed points of the
action. We can easily check that the
action is induced by a cyclic subgroup of
$\mathrm{SL}(2,\mathbb{C})$. Hence the singularities obtained are of type $A_{p-1}$. Together
this gives a result similar to the case of order 2. More precisely

\begin{prop}\label{p isogenie}
Let $S$ be a rational elliptic surface with semi-stable fibers and
a $p$ torsion section. Then there exists a $p$ isogeny between $S$
and some rational elliptic surface with section $S'$ such that:
\begin{itemize}
\item The surface $S'$ has only semi-stable fibers

\item If over some $t\in\mathbb{P}^1$ the fiber of $S$ is of type
$I_n$ then the fiber of $F'$ over $t$ is of type $I_{pn}$ or
$I_{\frac{n}{p}}$.
\end{itemize}
\end{prop}
\begin{rem}Observe that if there is an isogeny $\sigma:S_1\rightarrow
S_2$ then there is also an isogeny $\sigma:S_2\rightarrow S_1$.
\end{rem}
\subsection{Correspondences between fiber products}
Using above results we can construct correspondences between fiber
products in the following way.
\begin{cor}\label{corespondencje pomiedzy produktami}
Let $S_1$, $S_2$, $S'_1$ be rational elliptic surfaces with
section admitting only semi-stable fibers. Let $X$ and $X'$ be
small resolutions of $S_1\times_{\mathbb{P}^1}S_2$ and
$S'_1\times_{\mathbb{P}^1}S_2$. If $\sigma:S_1 \longrightarrow
S'_1$ is an isogeny introduced in Proposition \ref{p isogenie},
then
$(\sigma,\operatorname{id}):S_1\times_{\mathbb{P}^1}S_2\rightarrow
S'_1\times_{\mathbb{P}^1}S_2$ induces a correspondence between $X$
and $X'$.
\end{cor}
\begin{proof} It is clear that $(\sigma,\operatorname{id})$ restricts to a
finite map between open subsets of $X$ and $X'$.
\end{proof}
By composing correspondences of the above type we can construct a
large class of correspondences between Calabi--Yau three folds.

\begin{section}{Modular Calabi--Yau three folds}
Our aim is to prove the following.

\begin{conj} Let $S_1$ and $S_2$ be two non isogenous elliptic surfaces with
section defined over $\mathbb{Q}$ with semi-stable fibers such
that their fiber product $X=S_1 \times_{\mathbb{P}^1}S_2$ has
exactly five singular fibers. Let $\hat{X}$ be a small resolution
of $X$. Then there exists a correspondence between $X$ and some
rigid Calabi--Yau three fold $Y$. Moreover, $Y$ can be taken to be
the resolution of the double covering of $\mathbb{P}^3$ branched
over an octic surface.
\end{conj}
Unfortunately we will need some more assumptions. We now know two
methods that can lead to correspondences between nonrigid
Calabi--Yau fiber products and rigid Calabi--Yau three folds. The
first is described above, whereas the second was already
introduced in \cite{Schoen}. It is based on a fiberwise Kummer
construction.

\subsection{Correspondences with rigid fiber products. }
The first method gives correspondences between two fiber products.

Hence if the correspondence is defined over $\mathbb{Q}$ and one
of the three folds is rigid and satisfy the assumptions of the
theorem of Dieulefait and Maroharmayum then the modularity of both
varieties is proven.

To construct examples of correspondences we need first to analyze
isogenies between surfaces with at most 5 singular fibers. Let us
collect what is known about Beauville surfaces (i.e. surfaces with
exactly four singular fibers).

There are two possible birational descriptions of an elliptic
surface with semi-stable fibers as a double cover of
$\mathbb{P}^2$ ramified over a plane quartic.
\begin{itemize}
\item The first of them corresponds to the involution with respect
to the zero section. It leads to a quartic with a distinguished
point. The fibration is then given by inverse images of lines
passing through this point. This description is helpful in finding
explicitly 2-torsion sections.

\item The second corresponds to any involution of type $x\mapsto
b-x$ where $b$ is disjoint from the zero section. It leads to a
quartic and a point lying outside it. The fibration is then also
given by inverse images of lines passing through the point. This
description is convenient to find good Kummer fibrations.
\end{itemize}

The description of Beauville surfaces using quartics with
distinguished points is contained in \cite[table 6.7 and
6.8]{MirPer}.

The second description is in the following Table.
$$
\begin{tabular}{c|c|c}
\hline
name                    &fibers&   quartic with fibration (t:z) \\
\hline
 $\mathcal{E}(\Gamma(3))$ &3333    & $(x+t)(x^3-3tx^2+4z^3)$\\
 $\mathcal{E}(\Gamma_1(4)\cap\Gamma(2))$&4422& $(x+t+z)(x+t-z)(x-t+z)(x-t-z)$\\
$\mathcal{E}(\Gamma_1(5))$&5511&
$x(x^3-2x^2(z+t)+x(z^2+6zt+t^2)-4tz^2)$\\
$\mathcal{E}(\Gamma_1(6))$&6321& $(x-z)(x-t+2z)(x^2+t^2-z^2)$\\
$\mathcal{E}(\Gamma_0(8)\cap\Gamma_1(4))$&8211& $(x+z)(x-z)(x^2+t^2-z^2)$\\
$\mathcal{E}(\Gamma_0(9)\cap\Gamma_1(3))$&9111& $(x+z)(x^3-3tx^2+4z^3)$\\
\end{tabular}
$$

\begin{lem} \label{isogenie beau}Let $S$ be a Beauville surface with singular fibers
$I_{i_1}I_{i_2}I_{i_3}I_{i_4}$ over points $P_1,\dots,P_4\in
\mathbb{P}^1$. Then $S$ is isogenous to some Beauville surface
with a fiber of type $I_1$ over $P_1$.
\end{lem}
\begin{proof} We check all the possibilities. To find a torsion section
we use Lemma \ref{warunek na istnienie sekcji torsyjnej}. We
obtain the following.
\begin{itemize}
 \item The Beauville surface $\mathcal{E}(\Gamma(3))$ is isogenous
with the surface $\mathcal{E}(\Gamma_0(9)\cap \Gamma_1(3))$ as it
has a 3-torsion section. Moreover, analyzing more carefully we
observe that we also have an isogeny with
$\mathcal{E}(\Gamma_0(9)\cap \Gamma_1(3))^{\frac{1}{t}}$.
 \item The Beauville surface $\mathcal{E}(\Gamma_1(4)\cap\Gamma(2))$
 can be written after a suitable change of coordinates (without changing the
 fibration) as defined by the equation
 $$y^2=x(x-t)(xt-z^2).$$
 The fibration is given by the projection onto (z:t), hence a
 projection from the point (x=1,t=0,z=0), which is lying on the branch divisor of the
 involution $y\mapsto-y$. Observe that this situation is
 associated by a birational transformation to the double cover of
 $\mathbb{F}_2$ branched in the zero section and a 3-section.
 Moreover, the components of the 3-section correspond to components
 of the branch quartic. This means we have three different isogenies with three different surfaces.
 \item The Beauville surface $\mathcal{E}(\Gamma_1(5))$ admits a 5
 torsion section. We have only one possibility to take the
 quotient.
 \item The Beauville surface $\mathcal{E}(\Gamma_1(6))$ admits a 2
 torsion section and a 3 torsion section. Hence we have two
 possibilities to take the quotient
\end{itemize}
 If we moreover take into account that we can compose above isogenies we get the following
 table of groups of isogenous surface written in terms of
 corresponding fibers.
\begin{displaymath}
\begin{tabular}{c|c|c|c}
 3333& 4422&6231& 5511\\
 9111& 2244&2613& 1155\\
 1911& 8211&3162& \\
 1191& 2811&1326& \\
 1119& 1182&   & \\
     & 1128&   &
\end{tabular}
\end{displaymath}
\end{proof}

We can do something similar for surfaces with five singular
fibers. The following lemma is useful.
\begin{lem}[\cite{Looij}] \label{Looij} Elliptic surfaces with semi-stable
fibers are uniquely determined by their discriminant.
\end{lem}

After performing constructions as in the case of Beauville
surfaces we obtain.
\begin{lem}\label{isogenie 5 wlok} All elliptic surface with 5 singular fibers
can be divided into the following groups of isogenous surfaces
\end{lem}
\begin{displaymath}
\begin{tabular}{c|c|c|c|c|c}
 33321& 44211&62211& 54111& 53211&72111\\
 11163& 22422&31422&      &      &     \\
      & 11811&     &      &      &     \\
      & 11244&     &      &      &
\end{tabular}
\end{displaymath}
\begin{proof} The proof is analogous to the proof of Lemma \ref{isogenie
beau}. In the cases 42222 we need to find two two-torsion
sections. We use the fact that the sum $\sum \lfloor
\frac{k_i}{2}\rfloor =6$ i.e. we need to perform at least 6
blowings up to resolve the curve $T$. By the formula computing the
genus of a singular curve the number of components of $T$ is four.
\end{proof}
\begin{rem} Observe that in the case 62211 one of the fibers of
type $I_2$ is distinguished. This follows from the fact that a
quartic with an $A_5$ singularity cannot have three components.
Let us call the fiber that transforms to $I_1$ good and the one
that transforms into $I_4$ bad.
\end{rem}

Combining the two tables we obtain a large class of small
resolutions $\hat{X}$ of fiber products $X$ with exactly 5
singular fibers, that are isogenous to rigid Calabi--Yau fiber
products.

\begin{cor} \label{cor pomiedzy prod sztyw} Let $X$ be a fiber product of two elliptic surfaces $S_1$ and
$S_2$. If we can find representatives $S'_1$ and $S'_2$ of $S_1$
and $S_2$ in the tables of Lemmas \ref{isogenie beau},
\ref{isogenie 5 wlok} such that $S'_1 \times_{\mathbb{P}^1} S'_2$
has no fibers of type $I_0\times I_n$ for $n\geq2$, then $X$ is
isogenous to a rigid Calabi--Yau fiber product.
\end{cor}

\subsection{Correspondences with rigid Kummer fibrations.}
The second method is a fiberwise Kummer construction on fibers of
the fiber product. The obtained Kummer fibration is then a double
cover of $\mathbb{P}^3$ branched over an octic surface (for more
details see \cite{MG}).

Let us formulate a general context for this method. Let
$X=S_1\times_{\mathbb{P}^1}S_2$ be a fiber product of two elliptic
surfaces with section and with semi-stable fibers. Let $\hat{X}$
be its small projective resolution. Assume moreover that $X$ has
at most 5 singular fibers and no fibers of type $I_n\times I_0$
for $n\geq3$. We are interested in finding a rigid Calabi--Yau
Kummer fibration $\hat{Y}$ and a 2:1 rational map preserving
fibers $\phi:\hat{X}\dashrightarrow \hat{Y}$.

By \cite{CV} and \cite[Prop. 2.19.]{MG} we have two components of
the deformation space of any constructed Kummer fibration.
\begin{itemize}
\item The first is the so called transversal deformation space.
Its dimension depends on the sum of genera of the components of
the branch curve. It is zero in the case where all the components
of the branch curve of $\mathbb{\phi}$ are rational curves.

\item The second is the space of equisingular deformations. Its
dimension equals 0 if and only if the only fibers of type
$I_2\times I_0$ ($I_0\times I_2$) correspond to a node (not a
double tangent) of the quartic defining $S_1$ (resp. $S_2$).
\end{itemize}

Putting this together with Corollary \ref{cor pomiedzy prod sztyw}
we can prove the following theorem.
\begin{thm} \label{twierdzenie o przykladach corespondencji}
 Let $\hat{X}$ be a small resolution of
$X=S_1\times_{\mathbb{P}^1}S_2$ for some elliptic surfaces with
section $S_1$ and $S_2$. Assume one of the following:
\begin{itemize}
 \item The surfaces $S_1$ and $S_2$ are Beauville surfaces and $\sharp(S'')=3$.
 \item The surface $S_1$ is a Beauville surface, the surface $S_2$
 has five singular fibers, and $\sharp(S'')=4$. Moreover, $X$ does
 not admit a fiber of type $I_0\times I_5$ or $I_0\times I_7$, and
 does not admit a fiber of type $I_0\times I_6$ when $S_2$ is of
 type 62211.
\end{itemize}
Where $\sharp(S'')=3$ stands for the number of common singular
fibers of the surfaces $S_1$ and $S_2$. Then there is a
correspondence between $\hat{X}$ and some rigid Calabi--Yau three
fold $Y$.
\end{thm}
\begin{proof} By Corollary \ref{cor pomiedzy prod sztyw} we need only to consider the
case where $S_2$ has singular fibers $I_3I_3I_3I_2I_1$,
$I_4I_4I_2I_1I_1$ or $I_6I_2I_2I_1I_1$, and the product admits a
fiber of type $I_2\times I_0$ (in the last case we need to check
only the bad fiber). In these cases we construct rigid Kummer
fibrations. All of these surfaces can be seen as double covers of
$\mathbb{P}^2$ branched over quartic curves in such a way that the
fibers of type $I_2$ are induced by nodes. Let us be more precise.

We observe that each surface with fibers $I_3I_3I_3I_2I_1$
correspond to a quartic in $\mathbb{P}^2$ admitting three cusps
with a distinguished smooth point. By the Alekseev and Nikulin
description (cf. \cite{AN1} and Chapter 4) we can associate to
such a quartic the diagram of exceptional curves on the minimal
resolution of the double cover of $\mathbb{P}^2$ branched over the
quartic. In this case this diagram has to be the following.
$$ \setlength{\unitlength}{0.5mm}
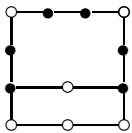
 \begin{picture}(70,30)(0,-20)
 \put(0,10){\circle{3}}
 \put(1.5,10){\line(1,0){7}}
 \put(8,8){$\bullet$}
 \put(18,8){$\bullet$}
 \put(30,10){\circle{3}}
 \put(11.5,10){\line(1,0){7}}
 \put(21.5,10){\line(1,0){7}}
 \put(0,8.5){\line(0,-1){7}}
 \put(0,-1.5){\line(0,-1){7}}
 \put(0,-11.5){\line(0,-1){7}}
 \put(-2,-2){$\bullet$}
 \put(-2,-12){$\bullet$}
 \put(28,-12){$\bullet$}
 \put(30,10){\circle{3}}
 \put(28,-2){$\bullet$}
 \put(30,-20){\circle{3}}
 \put(0,-20){\circle{3}}
 \put(15,-20){\circle{3}}
 \put(16.5,-20){\line(1,0){12}}
 \put(1.5,-20){\line(1,0){12}}
 \put(1.5,-10){\line(1,0){12}}
 \put(15,-10){\circle{3}}
 \put(16.5,-10){\line(1,0){12}}

 \put(30,8.5){\line(0,-1){7}}
 \put(30,-1.5){\line(0,-1){7}}
 \put(30,-11.5){\line(0,-1){7}}
\end{picture}
$$

The above graph gives us a way to construct the considered double
cover as a blow up of $\mathbb{P}^2$. We can hence associate to
the quartic a configuration of points (possibly infinitely near)
on $\mathbb{P}^2$. In this case the configuration of points is
uniquely determined by the position of four points which lie
generically on $\mathbb{P}^2$. This proves that the considered
quartic is unique up to an automorphism of $\mathbb{P}^2$. Hence
the surfaces with singular fibers $I_3I_3I_3I_2I_1$ form a
one-parameter family parameterized by the smooth points of the
quartic. We can construct a one parameter family of such surfaces
by taking double covers of $\mathbb{P}^2$ branched over a nodal
cubic and a line. By Lemma \ref{Looij} we construct in such a way
all surfaces with fibers of type $I_3I_3I_3I_2I_1$. The two
remaining cases can be dealt similarly. We obtain that all
rational elliptic surfaces with fibers $I_4I_4I_2I_1I_1$ or
$I_6I_2I_2I_1I_1$ are double covers of $\mathbb{P}^2$ branched
over a configuration of a conic and two lines.

Together with the result of \cite{CV} and \cite[Prop. 2.19.]{MG}
we proved that the corresponding Kummer fibrations have zero
dimensional equisingular deformations.

Let us compute their transversal deformations. The considered
quartics have respectively 2, 3 and 3 components. Using Lemma
\ref{isogenie beau} for each of the cases we need to check only
four possibilities corresponding to the four groups of Beauville
surfaces.

For economy of space we write explicitly only two examples, the
others being more or less analogous. For the surface with fibers
$I_4I_4I_2I_1I_1$ we consider the case corresponding to the third
column of the Table in the proof of Lemma \ref{isogenie beau}. We
can assume that the diagram of corresponding fibers is the
following.
$$\left(\begin{array}{ccccc}
4&4&2&1&1 \\
6&2& &3&1
\end{array}\right)
$$
In fact we should also consider the possibility where at the
bottom the 2 is exchanged with the 6, but it does not change
anything in what follows. In this case the involution restricted
to any fiber of the product has either 16 or 9 fixed points.
Moreover, the curve $C$ admits exactly two nodes. The Euler number
of the resolution of the curve $C$ is hence $-3 *16+3*16+18+2=20$.
The quartics describing both surfaces have 3 components. This
gives at least 9 components of the curve of intersection of the
two cones. We have at most 10 components as the intersection of
two quadric cones in which none is passing through the vertex of
the other can have at most two components. This implies that all
the components of the intersection of the quartic cones are
rational curves.

We consider also one example for the surface with fibers
$I_3I_3I_3I_2I_1$.

$$\left(\begin{array}{ccccc}
3&3&2&3&1 \\
8&2& &1&1
\end{array}\right)
$$
Here we have 16, 12 or 9 fixed points of the involution on each
fiber. The Euler number of the resolution of $C$ is 12. As the
intersection of the quartic cones has 6 components they all have
to be rational.

By checking all the cases we prove that the space of transversal
deformations is always zero-dimensional. Hence constructed Kummer
fibrations are rigid.
\end{proof}
\begin{rem}
It would be interesting to prove that none of the exceptions of
Theorem \ref{twierdzenie o przykladach corespondencji} is defined
over $\mathbb{Q}$.
\end{rem}
\end{section}

\vskip1cm
\begin{minipage}{7cm}
 Micha\l \ Kapustka\\
 Jagiellonian University\\
 ul. Reymonta 4\\
  30-059  Krak\'ow\\
 Poland\\
 email: Michal.Kapustka@im.uj.edu.pl
\end{minipage}
\end{document}